\documentclass[a4paper,12pt]{article}

\title{Quasi-invariant measures for some amenable groups  acting on the line}
\author{Nancy Guelman and Crist\'obal Rivas}
\date{}

\usepackage[top=2cm, bottom= 2cm, left=3cm, right=3.5cm, heightrounded,
  marginparwidth=3cm, marginparsep=3mm]{geometry}
\usepackage{diagbox}
\usepackage{psfrag}
\usepackage{fancyhdr}

\usepackage{tikz}
\usepackage{mathrsfs}
\usepackage{verbatim}
\usepackage[ansinew]{inputenc}
\usepackage{amsmath,amsthm,amssymb,amscd}
\usepackage{enumitem}

\newcommand{\N}{\mathbb{N}}
\newcommand{\Z}{\mathbb{Z}}

\newcommand{\R}{\mathbb{R}}

\newtheorem{rem}{Remark}

\newtheorem{teo}{Theorem}[section]

\newtheorem{prop}[teo]{Proposition}

\theoremstyle{definition}

\newtheorem{defi}{Definition}

\theoremstyle{remark}
\makeindex

\begin{document}

\maketitle

\begin{abstract}
In this note we show that if $G$ is a solvable group acting on the line, and if there is $T\in G$ having no fixed points, then there is a Radon measure $\mu$ on the line quasi-invariant under $G$. In fact, our method allows for the same conclusion for $G$ inside a class of groups that is closed under extensions and contains all solvable groups and all groups of subexponential growth.
\end{abstract}

%--------------------------------
\section{Introduction}
Let $G$ be a group acting by homeomorphism of the line.  We say that a (Borel) measure  $\mu$ is quasi-invariant or quasi-preserved under the action of $G$, if for every $g\in G$ there is $\lambda_g\in \R$ such that $g_*\mu=\lambda_g\mu$, where $g_*\mu: B\mapsto \mu(g^{-1}(B))$. We say that $\mu$ is preserved by $G$ if $g_*\mu=\mu$ for all $g\in G$.

Aiming to decide the (non-)amenability  of Thompson group $F$ (see \cite{thompson} for an introduction on this group, and \cite{juschenko} for an introduction to amenability) a very interesting criterion was proposed by L. Beklaryan  \cite[Theorem $B$]{beklaryan5}.

\vspace{0,3cm}

\noindent {\bf Criterion:} {\em If an amenable group acts by order-preserving homeomorphism of the line with an element acting freely, then  there is a Radon\footnote{A Borel measure $\mu$ on the line is said to be a Radon measure, if it gives finite mass to compact sets.} measure on the line quasi-invariant under the group action.}

\vspace{0,3cm}

Since for the natural -piecewise affine- action of $F$ on $(0,1)$ there is no quasi-invariant measure, the Criterion implies the non-amenability of $F$. However, the claim of validity of the Criterion was withdrawn \cite{beklaryan6},  apparently by the appearance of the preprint \cite{akhmedov} where it is claim that the Criterion fails already for the class of solvable groups.

This note is intended to clarify the discussion around the validity of the Criterion. We became interested in this problem after we discover a flaw in (the first version of) \cite{akhmedov}.  In fact, we will prove that the Criterion is valid in a class of groups that is closed under extensions and includes all solvable groups and all groups of subexponential growth (see \cite{de la harpe} for the definition of group growth).  We were, however, unable to decide weather Beklaryan criterion's holds in the class of all amenable groups.

We will say that a group $G$ has {\em locally} subexponential growth if any of its finitely generated subgroups has subexponential growth . Let $\mathcal{S}$ denote the class of groups $G$ for which there is a finite normal filtration
$$\{id\}=G^{d+1}\lhd G^d\lhd \ldots \lhd G^1\lhd G^0=G$$ with the property that $G_{i-1}/ G_i$ has locally subexponential growth, $i=1,\ldots, d$. Observe that $G_{i-1}/G_i$ may not be finitely generated. In this note, the {\em degree} of a group $G$ in $\mathcal S$ is the length of the shortest filtration in which each successive quotient has locally subexponential growth.  So for instance a group of subexponential growth has degree 1.

Clearly, any solvable group is in $\mathcal S$. Also, any group in $\mathcal S$ is amenable ($G$ in fact it is subexponentially elementary amenable, see \cite[Chapter 5]{juschenko}).  We will show
\vspace{0,3 cm}

\noindent{\bf Theorem A: }{\em Let $G$ be a group in $\mathcal{S}$ that is acting on the line by order-preserving homeomorphisms. Assume that there is $T\in G$ having no fixed points.
Then there is a Radon measure $\mu$ on the line which is quasi-preserved by $G$. }

\vspace{0,2 cm}

\begin{rem} In \cite{plante} Plante consider a class of groups $\mathcal{S}_0$ that contains all polycylic groups\footnote{A solvable group is polycylic if and only if it admits a filtration such that each successive quotient is cyclic.}, and all finitely generated groups of subexponential growth. He proves that {\em any} action on the line of a group in $\mathcal{S}_0$ quasi preserves a Radon measure. %In fact one can show that $\mathcal{S}_0$ contains all finite-rank solvable\footnote{A solvable group has  finite rank if each Abelian quitient appearing in its derived series embedds into $\Q^n$ for some $n$.} groups.

The class $\mathcal{S}_0$ however does not contains all solvable groups nor all groups that  have locally  subexponential growth. In fact counterexamples of Plante's theorem among finitely generated (infinite-rank) solvable group are easy to find: there are actions of some solvable groups not allowing for a quasi-invariant measure. For instance some actions of $\Z\wr \Z$ \cite{plante, rivas-tessera}, or some even more exotic as in \cite[\S 6.2]{botto-rhemtulla}. In these actions, though there are no global fixed points, each element of the group has at least one fixed point.

It is therefore natural to impose a priori in the Criterion, the condition that $G$ has an element acting without fixed points.
\end{rem}

Besides the groundwork provided by Plante, our main tool is the notion of {\em crossed elements}. This notion was introduce in \cite{beklaryan crossings}, but has been extensively studied/exploited in its connection with total orderings on groups (see \cite{navas order, rivas jgt, GOD}).
\begin{figure}[h] %h=here, t=top, b=bottom
\begin{center}
\begin{tikzpicture}

\draw[<->](-2.5,-2.5)->(2.5,2.5);
\draw[-] (-2,-2)-> (2.2,-0.5);
\draw[-] (-2.2,0.5)-> (2,2);

\draw (0.4,-1.5) node{$f$};
\draw (0,1.8) node{$g$};

\draw (2,1.6) node{$b$};
\draw (-2,-1.6) node{$a$};
\filldraw [black] (2,2) circle (1.3 pt);
\filldraw [black] (-2,-2) circle (1.3 pt);

\end{tikzpicture}
\end{center}
\caption{A crossing} \label{fig crossing}
\end{figure}
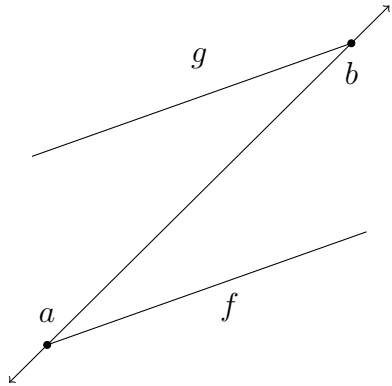
\begin{defi}We say $f$ and $g$,  two order-preserving homeomorphism of the line, are crossed if there is are $a,b\in\R$ such that $a<b$, $f(a)=a$, $g(b)=b$, and for every $x\in (a,b)$, $f(x)<x$ and $g(x)>x$. See figure \ref{fig crossing}, where the graphs of $f$ and $g$ are depicted.  
\end{defi}

For us, the main importance of crossed elements is that they entail exponential growth. Indeed, if $f$ and $g$ are crossed elements, then there is $n\in \N$ such that $f^n$, $g^n$ generates a free semi-group (see for instance \cite{navas}). In particular the group generated by $f$ and $g$ has exponential growth.

%\begin{enumerate}[label=(\roman*)]
%\item Since $A(G)$ quasi preserves the Lebesgue measure $Leb$, the group $G$ quasi-preserves the (Radon) measure $\rho^*Leb: B\mapsto Leb(\rho(B))$.
%\item Since $G$ (and hence $A(G)$) has no global fixed points,  $G$ {contains} an element having no fixed point.
%\end{enumerate}

%---------------------------------
\subsection{Quasi-invariant measures and semi conjugacy to affine actions}
\label{sec semiconj}
\begin{defi}Two representations $\rho_1,\rho_2:G\to Homeo_+(\R)$ are {\em semi conjugated} if there is a monotone map (i.e. non-decreasing) $c:\R\to \R$ which is proper (i.e. $c^{-1}$ sends compact sets to bounded sets, or, equivalently since $c$ is monotone, $c(\R)$ is unbounded in both directions) and such that for all $g\in G$
$$ c\circ \rho_1(g)=\rho_2(g)\circ c.$$
\end{defi}
\begin{rem} The above definition is the analog for actions of the line of the definition of semi conjugacy for groups acting on the circle from \cite{ghys lef}. Though, sometimes one also insists in the continuity of $c$ above (for instance in \cite{ghys,navas}), without the continuity assumption semi conjugacy becomes an equivalence relation (see \cite{ghys lef} for the case of the circle and \cite{ABR} for the case of the line). In this note, we do not assume continuity.
\end{rem}

Observe that for $G\subseteq Homeo_+(\R)$, the presence of a quasi-invariant measure $\mu$ provide us, for every $g\in G$, the affine map of the line
 \begin{equation}A_g(x)=\frac{1}{\lambda_g} x + \mu([0,g(0)),\label{eq afine}\end{equation}
where, by convention, if $g(0)<0$ then $\mu([0,g(0))$ is by definition $-\mu([g(0),0))$. This convention is made all across this note.

This association is in fact a representation of $G$ into the affine group:
\begin{eqnarray*}
A_{fg}(x)
&=& \frac{1}{\lambda_f}\left(\frac{1}{\lambda_g} x + \mu([f^{-1}(0),g(0))\right)= \frac{1}{\lambda_f}\left(\frac{1}{\lambda_g} x +\mu([0,g(0))+ \mu[f^{-1}(0),0)\right)\\
&=& \frac{1}{\lambda_f} \left(\frac{1}{\lambda_g}x+\mu([0,g(0))\right)+\mu([0,f(0))= A_f\circ A_g(x).
\end{eqnarray*}
Further, if the  $G$ action admits a quasi-invariant Radon measure $\mu$, then this action is semi conjugated to the affine action above. Indeed, if we let $F(x)=\mu([0,x))$, then
\begin{eqnarray*}
F (g(x))
&=& \mu([0,g(x)))\\
&=& \mu([g(0),g(x))+\mu([0,g(0))\\
&=& g^{-1}_*\mu([0,x))+\mu([0,g(0))\\
&=& \frac{1}{\lambda_g} F(x) + \mu([0,g(0))=A_g(F(x)).
\end{eqnarray*}
Observe that $\mu$ may have atoms, for instance when the $G$ action admits a discrete invariant set. We also have a converse

\begin{prop}If an action of $G$ is semi conjugated to an affine action, then there is a quasi-preserved measure for $G$. 
\end{prop}

\proof If $G$ has a global fixed point then the conclusion holds trivially, so we assume there are non. Suppose there is a semi conjugacy
$$F\circ g=A(g)\circ F,$$ for some affine representation $A:G\to Aff(\R)$. Since we are assuming that the $G$ action has no global fixed points, the same is true for the corresponding affine action of $G$.
Now, for a Borel set $B$, let $\mu(B):=Leb(F(B))$. We have that
\begin{eqnarray*}
g_*\mu(A)
&=& \mu(g^{-1}A)\\
&=& Leb(F\circ g^{-1}(A))\\
&=& Leb(A(g^{-1})\circ F(A))\\
%&=& \frac{1}{\lambda_g} Leb(F(A)) =\lambda_g\mu(A),
&=& \lambda_g Leb(F(A)) \;=\; \lambda_g\mu(A),
\end{eqnarray*}
where $\lambda_g$ is precisely the dilation factor of the map $A(g^{-1})$. $\hfill\square$

\section{Proof of Theorem A}

We begin with the next proposition which is the first step in an induction argument. The proposition is known, but we provide a full proof since the arguments in it will be use in the proof of Theorem A.
\begin{prop}
\label{prop d=1}

Let $G$ be a  subgroup of $Homeo_+(\R)$ locally of subexponential growth.  Assume that there is $T\in G$ having no fixed points.
Then there is a Radon measure $\mu$ on the line which is preserved by $G$.
\end{prop}

\noindent \textbf{Proof:} We first observe that the presence of the fixed-point-free element $T\in G$ implies that there is a non-empty minimal invariant set (that is, a closed set  invariant under $G$, and having no closed proper subset invariant under $G$).
Indeed, let $x_0$ be any point of the line and let $I=[x_0,T(x_0)], $ when $T(x_0)>x_0$ or  $I=[T(x_0),x_0], $ otherwise.
Let $\mathcal{F}$ the family of non-empty closed sets which are $G$-invariant.
Let us consider an order relation $\preceq$ in $\mathcal{F}$,  defined as $\Lambda_1 \succeq \Lambda_2$ if $\Lambda_1 \cap I  \subseteq \Lambda_2\cap I$.
Since $T$ has no fixed points, any $G$-orbit intersects $I$, therefore for any $\Lambda \in \mathcal{F}$, $\Lambda \cap I \neq \emptyset$.
By Zorn's lemma there exists a maximal element for $\preceq$. This maximal element is the intersection of $I$ with a minimal closed $G$-invariant set, which we denote by $\Lambda$. There are three possibilities.
\begin {itemize}
\item
$\Lambda '$ (the accumulation points of $\Lambda$) is empty.

In this case $\Lambda$ is discrete and $\mu=\sum_{m\in \Lambda} \delta_m$ is a Radon $G$-invariant measure.
\item
$\partial \Lambda$ (the boundary of $\Lambda$) is empty.

In this case, $ \Lambda= \R$, so the action of G is minimal (that is, every orbit is dense). We claim that this action is also free. Indeed, if the action is not free, then there exist $f \in G\setminus\{id\}$ having at least one fixed a point. We let  $I=[a,a')$, where $a\in \R$ and $a'\in\R\cup\{\infty\}$, be a non-empty component of $\R\setminus Fix(f)$. By eventually changing $f$ by its inverse, we can assume that $f(x)<x$, $\forall x \in I$.  Since the action is minimal, there is
$h \in  G$  such that $h(a) \in I$. Now consider the element $g=hf^n$. Since for any $x\in I$, $f^n(x)\to a$ as $n$ tends to $\infty$, we have that $g(x)>x$ for every $x\in [a,h(a)]$, but, if $n$ is large enough, we have that $g(f^{-1}h(a))=hf^{n} (f^{-1}h(a))<f^{-1}h(a)$. Thus $g$ has a fixed point that is greater than $h(a)$. Let $b$ be the infimum of these fixed points. Then $f$ and $g$ are crossed elements exactly as in Figure \ref{fig crossing}. This contradicts the fact that $G$ has locally subexponential growth.

%Let $J=[h(a),h(b))$.
%An argument analogous to last case of the proof of Theorem A shows that  $\langle g, h\rangle$ has exponential growth (in %fact $g$ and $h$ are crossed\marginpar{definir crossed} elements and by Lemma 2.2.44 of \cite{navas} $G$ contains a free %sub semi-group generated by two elements),and this contradicts that $G$ has subexponential growth.
%sending u or v inside ]u, v[; however, this implies that g and h are crossed on [u, v],
%thus contradicting our assumption.
Since the action of G is free and minimal, H\"{o}lder's theorem states that
 G is topologically conjugate to a  group of translations (see \cite{navas}).
We can pull back the Lebesgue measure by this conjugacy, to obtain an invariant Radon measure for the $G$-action.
\item $\Lambda'=\partial\Lambda=\Lambda.$
In this case  $\Lambda$  is  ``locally'' a Cantor set. So, one can collapses each interval in the complement of $\Lambda$ to a point to obtain another (topological) line and consider the induced $G$-action by {\em semi conjugacy}. 

More precisely, there is a surjective, non-decreasing and continuous map $c:\R\to \R$ that is constant in the complement of $\Lambda$. Since $\Lambda$ -and hence its complement- is $G$-invariant, we can define $\psi:G\to Homeo_+(\R)$ satisfying  
$$ \psi(g)\circ c= c\circ g\;,\;\; \forall g\in G.$$
In this way, since $c(\Lambda)=\R$, we have that $\psi(G)$ is group of homeomorphism acting minimally on the line: if $M$ is a $\psi(G)$-invariant set properly contained in $\R$, then $c^{-1}(M)\cap \Lambda$ is $G$-invariant properly contained in $\Lambda$. As, in the previous case, the group $\psi(G)$ also acts freely, so it preserves a Radon measure $\mu$. We can pull back this measure by the semi conjugacy by letting $\tilde{\mu}(A)=\mu(c(A))$, to obtain a $G$-invariant Radon measure.
\end{itemize}
$\hfill\square$

\noindent \textbf{Proof of Theorem A:} We argue by contradiction. Suppose there is a group $G$ in $\mathcal S$ contradicting the conclusion of Theorem A. We choose $G$ with the least possible degree. From Proposition \ref{prop d=1} the degree of $G$ is greater than 1, say it has degree $d+1$, and the filtration witnessing this degree is $\{id\}=G^{d+1}\lhd G^d\lhd \ldots \lhd G^1\lhd G^0=G$. We fix the action the action of $G$ on the line having no quasi-invariant measure, and also fix $T\in G$ an element having no fixed points. By eventually changing $T$ by its inverse, we can (and will) assume that $T(x)>x$ for all $x\in \R$. As in the proof of Proposition \ref{prop d=1} we have that $G$ has a minimal invariant set $\Lambda$.

If $\Lambda$ is discrete, then $\mu=\sum_{m\in \Lambda} \delta_m$ is a $G$-invariant measure, contradicting our assumption. So $\Lambda$ can not be discrete. We make two reductions.

\begin{enumerate}
\item We first argue that we can reduce to the case where  $G^d$ has no (global) fixed points.

Suppose $G^{d}$ has at least one fixed point. Since $G^{d}$ is normal in $G$, we have that $X:=\overline{Fix(G^{d})}$, the closure of set of $G^{d}$ fixed points, is an infinite $G$-invariant set unbounded in both directions of the line. Since $\Lambda $ is non-discrete, the closure of every $G$-orbit contains $\Lambda$  (see for instance \cite{ghys,navas}). Hence $\Lambda\subseteq X$.  

The action of $G$ on $X$ factor throughout an action of $G/G^{d}$. Moreover, the $G/G^d$ action on $X$ can be extended to an action on the whole real line, for instance by taking linear interpolation on the open components of the complement of $X$ %\footnote{In fact, the original action of $G$ is semi-conjugated (in the sense of Ghys) to the above action of $G/G^{d}$ on the line.}
(see for instance the proof Theorem 6.8 in \cite{ghys}). Denote this new action by $\psi:G\to Homeo_+(\R)$. Observe that in this construction $G^d$ acts trivially on the line and we have that $g(x)=\psi(g)(x)$ for every $x\in \Lambda$ and every $g\in G$. It follows that $\Lambda$ is also a minimal invariant set for $\psi(G)$.

From the minimality of the degree of $G$, we have that the action $\psi$ of $G/G^d$ on the line admits a quasi-invariant Radon measure $\mu$. In particular, $supp(\mu)$ the support of $\mu$, is closed and $G/G^d$-invariant, and hence $\Lambda\subseteq supp(\mu)$. We claim that $\Lambda=supp(\mu)$. To see this, first note that as in \S\ref{sec semiconj}, the presence of a quasi-invariant measure implies that $\psi(G)$ is semi conjugated to an affine group. Precisely, if we let $A_g$ as in (\ref{eq afine}) , and $F(x):=\mu([0,x))$, then
$$A_g\circ F= F\circ \psi(g)\, ,\; \forall g\in G.$$
This affine action has no global fixed points, and its orbits are not discrete (since $\Lambda$ is not discrete). It is therefore minimal: the closure of every orbit is the whole real line. Thus $\R=F(supp(\mu))=F(\Lambda)$. As a consequence we have that for any interval $I$ in the complement of $\Lambda$ there is $\alpha \in \R$ such that $F|_I=\alpha$. If we observe that two points $x$ and $y$ are identified under $F$ if and only if $\mu([x,y))=0$, we obtain that $I$ is  also in the complement of $supp(\mu)$. Therefore  $\Lambda=supp(\mu)$ as claimed.

The preceding claim implies, in particular, that $\Lambda=supp(\mu)\subseteq X$. Thus, in the original action of $G$,  $G^d$ fixes every point in $supp(\mu)$. Therefore $\mu$ is $G^d$-invariant, and hence $\mu$ is quasi-preserved by $G$. This contradict our choice of $G$. Hence, we conclude that $G^d$ has no global fixed points.

\item If there is an element  $T\in G^d$ acting freely, then Theorem A is ensured by Plante's \cite{plante} and Proposition \ref{prop d=1}.

Indeed, Proposition \ref{prop d=1} ensures the existence of a $G^d$-invariant Radon measure $\mu$ on the line. Since $T$ has no fixed points, the translation number homomorphism $\tau_\mu: g\mapsto \mu[0,g(0))$ defined on $G^d$ is non-trivial. Theorem A then follows from Lemma 4.1 and Lemma 4.2 in \cite{plante} (alternatively, the argument below also works).

\end{enumerate}

So we are left with the case where $G^d$ has no global fixed points but there is no element of $G^d$ acting freely. We claim that under the presence of the freely acting element $T\in G$ this case neither is possible.

Indeed, since $G^d$ has no global fixed points, there is $f\in G^d$ such that $T(0)<f(0)$. Since $f$ has at least one fixed point, there is a interval of the form $I=(a,b)$, where at least one of the endpoints is in $\R$ (and the other may be $\pm\infty$) containing $0$ fixed by $f$, and such that $f$ has no fixed point in its interior (hence $f(x)>x$ for all $x\in I$). For concreteness we assume that $a$ is a point in the real line, the other case being analogous.

Let $h=T^{-1}fT$. Since $G^d$ is normal in $G$, we have that $h\in G^d$. Observe that $h(a)\in I$. Then, proceeding in the same way that in the proof of Proposition \ref{prop d=1}, the case where $\partial \Lambda$  is empty, we can build $g\in G^d$ so that $f$ and $g$ are crossed. This contradicts that $G^d$ has locally subexponential growth. This last contradiction finishes the proof of Theorem A. $\hfill\square$

 %Let $J=[T(a),T(b)]$, it is easy to see that $J\cap I \neq \emptyset$, and for $h=Tg^{-1}T^{-1}$ $J$ is $h$-invariant (See figure 1).


\begin{thebibliography}{99}

\bibitem {ABR} {\sc J. Alonso, J. Brum, C. Rivas}. Flexibility of some subgroups of $Homeo_+(\R)$. Preprint available on arXiv:1605.07671.

\bibitem {botto-rhemtulla} {\sc R. Botto-Mura, A. Rhemtulla}. Orderable groups. Lecture Notes in Pure and Applied Mathematics, New York (1977).

\bibitem {akhmedov} {\sc A. Akhmedov}. Amenable subgroups of $Homeo(\R)$ with large characterizing quotient. http://arxiv.org/abs/1211.6165.

\bibitem {beklaryan5} {\sc L. Beklaryan}. The classification theorem for groups of homeomorphism of the line: the non-amenability of Thompson's groups $F$. http://arxiv.org/abs/1112.1942v5.

\bibitem {beklaryan6} {\sc L. Beklaryan}. The classification theorem for groups of homeomorphism of the line: the non-amenability of Thompson's groups $F$. http://arxiv.org/abs/1112.1942v6.


\bibitem {beklaryan crossings} {\sc L. Beklaryan}. Groups of homeomorphism of the line and the circle. Topological characteristics and metric invariants. {\em Uspehi Matem. Nauka} {\bf 59} (2004), 4-66. English translation: {\em Russian Math. Surveys} {\bf 59} (2004), 599-660.

\bibitem{thompson}{\sc  J. Cannon, W. Floyd, W. Parry.} Introductory notes on Richard Thompson's group. {\em L'einseignement Math\'ematique} {\bf 42} (1996), 215-256.

\bibitem{GOD} {\sc 
B. Deroin, A. Navas, C. Rivas.}
 {\em Groups Orders and Dynamics}. Preprint (2015 available on arXiv:1408.5805.

\bibitem{ghys lef}{\sc  E. Ghys.} Groupes d'homeomorphisms du cercle et cohomologie born�e. {\em Contemporary Mathematics} {\bf 58} (1987), part III, 81-106. 

\bibitem{ghys}{\sc  E. Ghys.} Groups acting on the circle. {\em L'einseignement Math\'ematique} {\bf 42} (2001), 329-407.

\bibitem{de la harpe} {\sc P. de la Harpe} .
 {\em Topics in geometric group theory}. Chicago Lectures in Mathematics.


\bibitem{juschenko} {\sc K. Juschenko.}
 {\em Amenability}. In preparation (2015). Current version available at 
 {\tt http://www.math.northwestern.edu/\string~juschenk/book.html}.

\bibitem {navas} {\sc A. Navas.} {\em Groups of circle diffeomorphism}. Chicago Lectures in Mathematics.

\bibitem {navas order} {\sc A. Navas.} On the dynamics of (left) orderable groups. {\em Ann. Inst. Fourier (Grenoble)} {\bf 60} (2010), 1685-1740.

\bibitem {plante} {\sc J.F. Plante.} Solvable groups acting on the line, Trans. Amer. Math. Soc (1983), 401-414.

\bibitem {rivas jgt} {\sc C. Rivas.} On the space of Cornadian group orderings. {\em J. Group Theory} {\bf 13} (2010), 337-353.

\bibitem {rivas-tessera} {\sc C. Rivas, R. Tessera.} On the space of left-orderings of virtually solvable groups. {\em Groups, Geometry and Dynamics} {\bf 10} (2016), 65-90.


\end{thebibliography}
\end{document}